\documentclass[12pt, reqno]{amsart}
\usepackage{amsmath, amsthm, amscd, amsfonts, amssymb, graphicx, color}
\usepackage[bookmarksnumbered, colorlinks, plainpages]{hyperref}
\newtheorem{thm}{Theoreme}[section]

 \theoremstyle{definition}
 
 \theoremstyle{remark}

 \numberwithin{equation}{section}
\begin{document}
\title[Biharmonic Immersion]{Biharmonic Immersion in Cartan–Hadamard}
\author{Sa\"{\i}d Asserda and M'Hamed Kassi}
\address{Ibn tofail university , faculty of sciences, department of mathematics, PO 242 Kenitra Morroco}
\email{asserda-said@univ-ibntofail.ac.ma}

\address{Regional Centre of trades of education and training, Kenitra Morocco}
\email{M'hamedkassi@yahoo.fr}
\date{\today}
\begin{abstract}
If $(N^{m+p},h)$ is a Cartan-Hadamard manifold such that $Ric(h)\geq -G(r_{N}(x))$ where $G(0)\geq 1, G^{'}\geq 0$ and $G^{-1/2}\not\in L^{1}(+\infty)$ then every  proper biharmonic isometric immersion  $\phi : M^{m}\rightarrow(N^{m+p},h)$  is a harmonic map.
\end{abstract}
\maketitle
\section{Introduction}
Soit $ \phi : (M,g)\rightarrow (N,h)$ de classe $C^{\infty}$ d'une vari\'et\'e Riemannienne $M$ de dimension $m$ dans $N$ de dimension $n$. L'\'energie   et la bi-en\'ergie de $\phi$ sur $\Omega\subset\subset M$ sont d\'efinie par
$$
E_{1}(\phi)=\int_{\Omega}<d\phi,d\phi>_{h}dv_{g},\qquad E_{2}(\phi)=\int_{\Omega}<\tau(\phi),\tau(\phi)>_{h}dv_{g}
$$
o\`u $\tau(\phi)=\hbox{trace}_ {g}\nabla d\phi$ est le champ de tension de $\phi$ et $dv_{g}$ est la forme volume de $M$. L'application $\phi$ est harmonique ( minimale  ) si $\phi$ est une extr\'emale de $E_{1}$ i.e $\tau(\phi)=0$ et biharmonique ( o\`u biminimal ) si $\phi$ est une extr\'emale de  de $E_{2}$ i.e  $\tau_{2}(\phi):=-\Delta^{\phi}\tau(\phi)-\hbox{trace}_{g}R^{N}(d\phi,\tau(\phi))d\phi=0$ o\`u $\Delta^{\phi}=-\hbox{trace}_{g}(\nabla^{\phi}\nabla^{\phi}-\nabla^{\phi}_{\nabla^{M}})$ est le laplacien sur les sections du fibr\'e $\phi^{-1}(TN)$ au dessus de $M$  et $R^{N}$ est le tenseur courbure de $(N,h)$ d\'efini par $R^{N}(X,Y)Z:=\nabla_{X}\nabla_{Y}Z-\nabla_{Y}\nabla_{X}Z-\nabla_{[X,Y]}Z$. Une sous vari\'et\'e $M\subset N$ est minimale ( resp. biminimal ) si l'injection de $M$ dans $N$ l'est. Dans [1], B.Y.Chen a pos\'e la conjecture : toute sous vari\'et\'e bi-minimale de l'espace euclidien est minimale. Il expose dans  son  survey les d\'eveloppements r\'ecents sur la conjecture [2]. Dans [4], S.Maeta   propose une version g\'eom\'etrique globale de la conjecture de Chen : toute immersion compl\`ete bi-minimale dans une vari\'et\'e riemannienne de courbure sectionnelle  n\'egative  est minimale. Il montre que c'est le cas si l'immersion est propre ( en particulier compl\`ete ) et la courbure sectionnelle  $K_{N}$ v\'erifie $-C(1+d_{N}^{2}(x,x_{0}))^{\alpha\over 2}\leq K_{N}(x)\leq 0$ o\`u $C\geq 0$ et  $0\leq\alpha<2$.\\
Dans cette note,  en supposant que l'espace ambiant $(N,h)$ est de Cartan-Hadamard et sa  courbure de Ricci  d\'ecroit  vers $-\infty$ au plus quadratique en la distance, on montre que l'immersion est minimale. 
\begin{thm}
Soit $\phi : M^{m}\rightarrow (N^{m+p},h)$ une immersion isom\'etrique  propre o\`u $N$ une vari\'et\'e de Cartan-Hadamard de courbure de Ricci :  $Ric_{h}(x)\geq -G(d_{N}(x,x_{0}))$ avec $G : [0,+\infty[\rightarrow[0,+\infty[$ v\'erifiant  $G(0)\geq 1,\ G^{'}\geq 0$ et $G^{-1/2}\not\in L^{1}([0,+\infty[)$. Si $\phi$ est bi-minimale alors elle est minimale.
\end{thm}
\section{Pr\'eliminares}
Sot $ M^{m}\subset(N^{m+p},h)$ une sous vari\'et\'e  de dimension $m$ dans une vari\'et\'e de dmension $m+p$. On munit $M$ de la m\'etrique induite par l'injection canonique $i_{M}$ et $h$ : $i_{M}^{*}(h)(X,Y):=<di_{M}(X),di_{M}(Y)>_{h}$ pour $X,Y\in TM$. La deuxi\`eme forme fondamentale de $M$ : $ B : TM\times TM\rightarrow NM$ est d\'efinie par
$$
B(X,Y)=D_{X}Y-\nabla_{X}Y
$$
pour tout $XY \in TM$ o\`u $D$ est la connexion de Levi-Civita de $(N,h)$, $\nabla$ est celle de $(M,i_{M}^{*}h)$ et $NM=TN\circleddash TM$ est le fibr\'e normale de $M$. Si $\eta\in TN$, l'application de Weingarten associ\'ee a $\eta$ : $A_{\eta} : TM\rightarrow TM$ est d\'efinie par $D_{X}\eta=A_{\eta}X+\nabla^{\perp}_{X}\eta$
o\`u $\nabla^{\perp}$ est la connexion normale, elle v\'erifie   $$<B(X,Y),\eta>_{h}=<A_{\eta}X,Y>_{h}$$ Si $x\in M$, soit $(e_{1},e_{2},\cdots,e_{m},e_{m+1},\cdots,e_{m+p})$  une base locale orthonormale de $T_{x}N$ telle que $(e_{1},e_{2},\cdots,e_{m})$ est une base orthonormale de $T_{x}M$.  Alors $B$ se d\'ecompose au point $x$ : $B(X,Y)=\displaystyle\sum_{\alpha=m+1}^{m+p}B_{\alpha}(X,Y)e_{\alpha}$.
La courbure moyenne de $M$ est la trace de $B$
$$
H={1\over m}\sum_{i=1}^{m}B(e_{i},e_{i})=\sum_{\alpha=m+1}^{m+p}H_{\alpha}e_{\alpha}\ \
\hbox{o\`u}\ \ H_{\alpha}={1\over m}\sum_{=1}^{m}B_{\alpha}(e_{i},e_{i})
$$
Maintenant, soit $\phi : M\rightarrow(N^{m+p},h)$ est une immersion isom\'etrique. On identifie $d\phi(X)$ avec $X\in T_{x}M$ pour tout $x\in M$. Pour tous $X,Y \in TM$, la deuxi\`eme forme fondamentale de $\phi$ est $\nabla d\phi(X,Y)=\nabla^{\phi}_{X}(d\phi(Y))-d\phi(\nabla_{X}Y)=B(X,Y)$
et par suite $\tau(\phi)=mH$.
Donc $\phi$ est biharmonique si et seulement si $-\Delta H-\sum_{i=1}^{m}R^{N}(e_{i},H)e_{i}=0$, qu'on  d\'ecompose en parties tangentielle et normale ( voir [1]) :
$$
(1)\quad\qquad\Delta^{\perp}H-\sum_{i=1}^{m}B(A_{H}e_{i},e_{i})+\sum_{i=1}^{m}(R^{N}(e_{i},H)e_{i})^{\perp}=0
$$
$$
(2)\quad\qquad m\nabla\vert H\vert^{2}+4\sum_{i=1}^{m}A_{\nabla^{\perp}_{e_{i}}H}e_{i}-\sum_{i=1}^{m}(R^{N}(e_{i},H)e_{i})^{\top}=0
$$
\section{D\'emonstration du th\'eor\`eme 1.1}
\noindent  Puisque $\phi$ est bi-minimale, d'apr\`es l'equation $(1)$ on a
$$\Delta\vert H\vert^{2}=2\vert\nabla^{\perp}H\vert^{2}+2<H,\Delta^{\perp}H>\qquad\qquad\quad\qquad\qquad\qquad$$
$$\qquad\qquad\qquad\quad\qquad=2\vert\nabla^{\perp}H\vert^{2}+2\sum_{i=1}^{m}<B(A_{H}e_{i},e_{i}),H>-2\sum_{i=1}^{m}<R^{N}(e_{i},H)e_{i},H>$$
$$\qquad\geq2\vert\nabla^{\perp}H\vert^{2}+2\sum_{i=1}^{m}<A_{H}e_{i},A_{H}e_{i}>\quad\hbox{car $R^{N}\leq 0$}
$$
$$= 2\vert\nabla^{\perp}H\vert^{2}+2m\vert H\vert^{4}\qquad\qquad\qquad\qquad\qquad\quad
$$
On pose $u(x)=\vert H(\phi(x))\vert^{2}$ et on consid\`ere la fonction
$$
F(x)=(R^{2}-r^{2}(\phi(x)))^{2}u(x)\quad\hbox{si}\quad x\in M\cap\phi^{-1}(\overline{B}_{R})
$$
o\`u $y_{0}\in N\setminus\overline{\phi(M)}$ est un point fix\'e, $r(\phi(x))=d_{h}(\phi(x),y_{0})$ et $B_{R}=\{ y\in N\ :\ d_{h}(y,y_{0})<R\}$. La fonction $F$ est non identiquement nulle sur $M\cap\phi^{-1}(\overline{B}_{R})$ et nulle sur $M\cap\phi^{-1}(\partial\overline{B}_{R})$. Puisque $\phi$ est propre $M\cap\phi^{-1}(\overline{B}_{R})$ est compacte dans $M$. Il existe un maximum $p\in M\cap\phi^{-1}(B_{R})$ de $F$. Par un argument de Calabi, on peut supposer que $\phi(p)$ n'est pas conjugu\'e \`a $y_{0}$. On a donc
$$
(i)\qquad \nabla F(p)=0\Longleftrightarrow{\nabla F(p)\over F(p)}={2\nabla r^{2}(\phi(p))\over R^{2}-r^{2}(\phi(p))}\quad\hbox{et}\quad(ii)\ \ \Delta F(p)\leq 0
$$
et en utilisant $(i)$, $(ii)$ devient
$$
(iii)\qquad{\Delta u(p)\over u(p)}\leq{6\vert\nabla r^{2}(\phi(p))\vert^{2}\over(R^{2}-r^{2}(\phi(p)))^{2}}+{2\Delta r^{2}(\phi(p))\over R^{2}-r^{2}(\phi(p))}
$$
On a $\vert\nabla r^{2}(\phi(p))\vert^{2}\leq 4mr^{2}(\phi(p))$ et
$$
\Delta r^{2}(\phi(p))=2\sum_{=1}^{m}<(\nabla r)(\phi(p)),d\phi(e_{i})>^{2}+2r(\phi(p))\sum_{i=1}^{m}\nabla dr(\phi(p))(d\phi(e_{i}),d\phi(e_{i}))
$$
$$
+2r(\phi(p))<(\nabla r)(\phi(p)),\tau(\phi)(\phi(p))>
$$
$$ \leq 2m+2r(\phi(p))\Bigl(\sum_{i=1}^{m+p}\nabla dr(\phi(p))(e_{i},e_{i})-\sum_{i=m+1}^{m+p}\nabla dr(\phi(p))(e_{i},e_{i})\Bigr)$$
$$+2mr(\phi(p))\vert H(\phi(p)\vert
$$
Puisque $N$ est de Cartan-Hadamard on a $\nabla dr(\phi(p))(e_{i},e_{i})\geq 0$ (voir [3]),  par suite
$$
\Delta r^{2}(\phi(p))\leq 2m+2r(\phi(p))\Delta r(\phi(p))+2mr(\phi(p))\vert H(\phi(p)\vert
$$
D'apr\`es le th\'eor\`eme de comparaison : $\Delta r(\phi(p))\leq C\sqrt{G(r(\phi(p))}$ (voir  [5]), donc $(iii)$ s'\'ecrit
$$
{\Delta u(p)\over u(p)}\leq {24mr^{2}(\phi(p))\over(R^{2}-r^{2}(\phi(p)))^{2}}+{2m+2Cr(\phi(p))\sqrt{G(r(\phi(p)))}+2mr(\phi(p))\vert H(\phi(p)\vert\over  R^{2}-r^{2}(\phi(p))}
$$
D'apr\`es $(*)$ on a $\Delta u\geq 2mu^{2}$ sur $M$, on a donc
$$
2mu(p)\leq {24mr^{2}(\phi(p))\over(R^{2}-r^{2}(\phi(p)))^{2}}+{2m+2Cr(\phi(p))\sqrt{G(r(\phi(p)))}+2mr(\phi(p))\vert H(\phi(p)\vert\over  R^{2}-r^{2}(\phi(p))}
$$
et par suite
$$
2mF(p)\leq 24mr^{2}(\phi(p))+(2m+2Cr(\phi(p))\sqrt{G(r(\phi(p)))})(R^{2}-r^{2}(\phi(p))+2mr(\phi(p))\sqrt{F(p)}
$$
Puisque $r(\phi(p))\leq R$ et $G\ge 1$, cette in\'egalit\'e quadratique en $\sqrt{F(p)}$ entraine
$$
\vert H(\phi(x))\vert\leq C\sqrt{G(R)}\qquad\forall\ x\in M\cap\phi^{-1}(\overline{B_{R}}),\quad \forall\ R > 0.
$$
Si $w\in M$ et $R=d_{h}(\phi(w),y_{0})$, alors on a
$$
\vert H(\phi(w))\vert\leq C\sqrt{G(r(\phi(w)))}\qquad\forall\ w\in M
$$
et donc
$$
\Delta(r\circ\phi)^{2}\leq C(r\circ\phi)\sqrt{G(r\circ\phi)}\qquad\hbox{sur}\quad M
$$
On en d\'eduit que $(G,(r\circ\phi)^{2})$ est une paire d'Omori-Yau dans $M$, par suite $\Delta_{M}$ v\'erifie le principe de maximum [5].\\
Supposons que $u$ est major\'e sur $M$, pour tout $\epsilon > 0$ il existe $x_{\epsilon}\in M$ tels que
$$
\sup_{M}u\leq u(x_{\epsilon})+\epsilon,\quad\vert\nabla u(x_{\epsilon})\vert < \epsilon\quad\hbox{et}\quad\Delta u(x_{\epsilon})<\epsilon
$$
D'apr\`es $(*)$ on a $2mu^{2}\leq\Delta u$, donc
$$
0\leq\sup_{M}u\leq u(x_{\epsilon})+\epsilon\leq ({1\over 2m}+1)\epsilon
$$
par suite $u=\vert H\vert^{2}=0$ i.e $M$ est minimale.\\
Supposons que $u$ n'est pas major\'e. La fonction  $v=(u+1)^{-1/4}$ est minor\'ee sur $M$. D'apr\`es le principe de maximum : pour tout $n \geq 1,\ \exists\ x_{n}\in M$ tels que
$$
v(x_{n}) < \inf_{M}v+{1\over n},\ \ \vert\nabla v(x_{n})\vert < {1\over n},\ \ \Delta v(x_{n}) > -{1\over n}
$$
Puisque
$$
(1+u)^{-3/2}\Delta u=-4(1+u)^{-1/4}\Delta v + 20\vert\nabla g\vert^{2}
$$
et $2m^{2}u\leq\Delta u$, on en d\'eduit
$$
{2m(u(x_{n}))^{2}\over (1+u(x_{n}))^{3\over 2}} < {4\over n(1+u(x_{n}))^{1\over 4}}+{20\over n^{2}}
$$
Quand $n$ tend vers l'infini, $u(x_{n})\rightarrow\sup u=+\infty$. Le membre de droite tend vers l'infini et celui de gauche vers $0$, impossible.

\end{document}